\documentclass[11pt]{amsart2000}
\newtheorem{theorem}{Theorem}[section]
\theoremstyle{definition}
\newtheorem*{example}{Example}

\begin{document}
\setlength{\unitlength}{0.01in}
\linethickness{0.01in}
\begin{center}
\begin{picture}(474,66)(0,0) 
\multiput(0,66)(1,0){40}{\line(0,-1){24}}
\multiput(43,65)(1,-1){24}{\line(0,-1){40}}
\multiput(1,39)(1,-1){40}{\line(1,0){24}}
\multiput(70,2)(1,1){24}{\line(0,1){40}}
\multiput(72,0)(1,1){24}{\line(1,0){40}}
\multiput(97,66)(1,0){40}{\line(0,-1){40}} 
\put(143,66){\makebox(0,0)[tl]{\footnotesize Proceedings of the Ninth Prague Topological Symposium}}
\put(143,50){\makebox(0,0)[tl]{\footnotesize Contributed papers from the symposium held in}}
\put(143,34){\makebox(0,0)[tl]{\footnotesize Prague, Czech Republic, August 19--25, 2001}}
\end{picture}
\end{center}
\vspace{0.25in}
\setcounter{page}{309}
\title[lsc multifunctions in quasi-uniform and vector spaces]{On lower 
semicontinuous multifunctions in quasi-uniform and vector spaces}
\author{Andrzej Spakowski}
\address{Instytut Matematyki i Informatyki, Uniwersytet Opolski,\\
Oleska 48, 45-052 Opole, Poland}
\email{aspakow@uni.opole.pl}
\thanks{Andrzej Spakowski,
{\em On lower semicontinuous multifunctions in quasi-uniform and vector 
spaces},
Proceedings of the Ninth Prague Topological Symposium, (Prague, 2001),
pp.~309--319, Topology Atlas, Toronto, 2002}
\begin{abstract}
Given a cover $\mathcal{B}$ of a quasi-uniform space $Y$ we introduce 
a concept of lower semicontinuity for multifunctions $F:X\to 2^Y$, 
called $\mathcal{B}$-lsc. In this way, we get a common
description of Vietoris-lsc, Hausdorff-lsc, and bounded-Hausdorff-lsc 
as well. 
Further, we examine set-theoretical and vector operations on such 
multifunctions. We also point out that the convex hull of Hausdorff-lsc 
multifunctions need not to be Hausdorff-lsc except the case where the 
range space is locally convex.
\end{abstract}
\maketitle

\section{Lower semicontinuous multifunctions}

The two most known concepts of lower semicontinuity for multifunctions 
are the lower semicontinuity in Vietoris sense (V-lsc) and the lower 
semicontinuity in Hausdorff sense (H-lsc).
Given a set $Y$ we denote by $2^Y$ the family of all subsets
of $Y$. Every map $F:X\rightarrow 2^Y$ will be called
a {\it multifunction} from $X$ to $Y$.
Now, let $X$ and $Y$ be two arbitrary topological spaces.
We say that a multifunction $F:X \rightarrow 2^Y$ is {\it V-lsc}
at a point $x_0 \in X$ provided for every open $G \subset Y$
such that $F(x_0) \cap G \neq \emptyset $ there exists a neighbourhood
$U(x_0)$ of $x_0$ such that
$F(x)\cap G\neq\emptyset$ for every $x\in U(x_0)$.
This is the first concept of lower semicontinuity.

Let $(Y,{\mathcal{U}})$ be a uniform space. Recall that every uniformity
generate a topology, and a topological space is uniformizable
provided it is a Tichonov space. For multifunctions from $X$ to
$(Y,\mathcal{U})$ we may formulate the second concept of lower
semicontinuity. Namely, a multifunction $F:X \rightarrow 2^Y$ is called
{\it H-lsc at} $x_0 \in X$ if for every $W\in{\mathcal{U}}$ there exists
a neighbourhood $U(x_0)$ of $x_0$ such that
$$
F(x_0) \subset W(F(x))\ \text{for every}\ x\in U(x_0),
$$
where 
$$W(F(x))=\{\ y \in Y: (z,y) \in W \textrm{ for some } z \in F(x)\}.$$
This is the second concept of lower semicontinuity.
In particular, if $Y$ is a topological vector space, with its natural
uniformity generated by the neighbourhoods of 0, the condition of H-lsc
can be written in the equivalent form:
$$
F(x_0) \subset F(x)+V\ \text{for every}\ x\in U(x_0),$$
where V is a neighbourhood of 0 and 
$C+D=\{ c+d:c \in C,\, d \in D\}$
is the vector sum of sets $C$ and $D$.

It is known that every topological space $Y$ is
{\it quasi-uniformizable} (\cite{Pervin-1962, Murdeshwar-Naimpally-1966}).
This means that there is family ${\mathcal{U}}$ of subsets of
$Y\times Y$ such that:
\begin{enumerate}
\item
every $U\in{\mathcal{U}}$ contains the diagonal $\Delta$ of $Y\times Y$,
\item
$U,V \in {\mathcal{U}}$ implies that $U \cap V \in {\mathcal{U}}$,
\item
for every $U \in {\mathcal{U}}$ there exists $V \in\mathcal{U}$ such that
$V \circ V \subset U$, where $V \circ V = \{ (x,y) \in Y \times Y:(x,z),
(z,y) \in V \textrm{ for some } z\in Y \}$,
\item
$U \in\mathcal{U}$ and $U \subset V$ implies $V \in\mathcal{U}$,
\item
the family $\{ W(y):W \in\mathcal{U},\ y \in Y \} $ is a neighbourhood 
system generating the topology of $Y$, where 
$W(y)= \{ z \in Y:(y,z) \in W \}$.
\end{enumerate}
Every such family is called a {\it quasi-uniformity of the topological
space} $Y$.

If $Y$ is only a set, then a family $\mathcal{U}$ of subsets of $Y\times Y$
satisfying properties (1)--(4) is called a {\it quasi-uniformity on} $Y$,
and the pair $(Y,\mathcal{U})$ a {\it quasi-uniform space}.
A quasi-uniform space $(Y,\mathcal{U})$ is a uniform space provided
$\mathcal{U}$ have the following symmetric property:
$W\in\mathcal{U}$ implies $W^{-1}\in\mathcal{U}$,
where $W^{-1}=\{ (z,y)\in Y\times Y:(y,z)\in W \}$.

In every topological space $Y$ we have:
for every $A\subset Y$ and every quasi-uniformity $\mathcal{U}$
of $Y$,
$\textrm{\it cl\,}(A)=\bigcap\{ W^{-1}(A): W\in\mathcal{U} \}$.
where $\textrm{\it cl\,}(A)$ denotes the closure of $A$ in $Y$.
The usage of $W^{-1}$ is explained by the following facts:
\begin{enumerate}
\item
the sets $W(y)$, where $W\in\mathcal{U}$ and $y\in Y$, form a
neighbourhood system for the topology of $Y$,
\item
for arbitrary $A\subset Y$ we have: $y\in\textrm{\it cl\,}(A)$, the
closure of $A$, if and only if for every $W\in\mathcal{U}$, $W(y)\cap
A\neq\emptyset$ if and only if there exists $a\in A$ such that $a\in
W(y)$, or equivalently $y\in W^{-1}(a)$.
\end{enumerate}
For a quasi-uniform $(Y,\mathcal{U})$ space the definition of H-lsc
should be modified as follows. 
A multifunction $F:X \rightarrow 2^Y$ is called {\it H-lsc at} $x_0 \in X$
if for every $W \in \mathcal{U}$ there exists a neighbourhood $V$ of $x_0$
such that 
$$
F(x_0) \subset W^{-1}(F(x))\ \text{for every}\ x\in V.
$$

Note that we have the following property: $F$ is H-lsc at $x_0$
if and only if $\textrm{\it cl\,}(F)$ is H-lsc, where $\textrm{\it cl\,}(F)$
is the closure multifunction of $F$,
i.e.\ $\textrm{\it cl\,}(F)(x)=\textrm{\it cl\,}(F(x))$ for all $x\in X$.
The basic relationships between V-lsc and H-lsc are well-known
(see 
\cite{Klein-Thompson-1984},
\cite{Lechicki-1979}).
Namely, if $Y$ is a topological space, $\mathcal{U}$ a quasi-uniformity
of $Y$ and $F:X \rightarrow 2^Y$ a multifunction, then H-lsc of
$F$ at $x_0\in X$ implies its V-lsc at $x_0$.
The converse holds provided the set $F(x_0)$ is totally bounded.
Recall that $A\subset Y$ is called {\it totally bounded}
provided for every $W\in \mathcal{U}$ there exists a finite set
$B\subset A$ such that $A\subset W^{-1}(B)$.
In general, V-lsc need not imply H-lsc.

\section{$\mathcal{B}$-lower semicontinuity}

Penot \cite{Penot-1993} introduced a concept of bounded lower semicontinuity
for multifunction from a topological space to a normed space.
In \cite{Borwein-Vanderwerff-1996} a similar idea is applied to convergence 
of sets, in particular, to convergence of epigraphs, 
with respect to the families of single subsets, compact subsets, 
weakly compact subsets, and of bounded subsets.
Following this, we define an abstract concept of lower semicontinuity
to unify the description of the above mentioned semicontinuities.

Let $X$ be a topological space, $(Y,\mathcal{U})$ a quasi-uniform
space, $\mathcal{B}$ a {\it cover} of $Y$, i.e., a family of nonempty
subsets of $Y$ such that $\bigcup\mathcal{B}=Y$.
We say that $F:X\rightarrow 2^Y$ is {\it $\mathcal{B}$-lsc}
at $x_0\in X$ provided for every $W\in \mathcal{U}$ and every
$B\in\mathcal{B}$
there exists a neighbourhood $U(x_0)$ of $x_0$ such that
$$
F(x_0)\cap B\subset W^{-1}(F(x))\ \text{for every}\ x\in U(x_0).$$

If $F(x_0)=\emptyset$, then $F$ is trivially $\mathcal{B}$-lsc at $x_0$ for 
arbitrary cover $\mathcal{B}$. Note also the following simple observations
and remarks:
\begin{enumerate}
\item
If $\mathcal{B}=\{ Y \}$, then $\mathcal{B}$-lsc is simply the
H-lsc.
\item
If $\mathcal{B}$ is a cover of $Y$ and
$\mathcal{B}_1$ the family of all finite unions of subsets of
$\mathcal{B}$, then $\mathcal{B}$-lsc implies the $\mathcal{B}_1$-lsc.
\item
If $\mathcal{B}\subset\mathcal{B}_1$, then $\mathcal{B}_1$-lsc
implies $\mathcal{B}$-lsc.
\item
If $\mathcal{B}$ is the family of all balls $B(r)=\{ y\in Y:||y||<r\}$
of a normed space $Y$, then $\mathcal{B}$-lsc will be called also 
{\it bounded H-lsc}. This case is identical with the Penot's 
concept \cite{Penot-1993}.
\item
If for every $B\in\mathcal{B}$ the multifunction $F_B(x)=F(x)\cap B$, 
$x\in X$, is H-lsc at $x_0$, then it is $\mathcal{B}$-lsc at this point. 
The converse, in general, does not hold. For instance, we can take:
$F(0)=[0,1]$ and $F(x)=[0,1)$ for $x\in (0,1]$, and $\mathcal{B}$
the family of all singletons of $[0,1]$. Then $F$ is $\mathcal{B}$-lsc but 
$F_B$, where $B=\{ 1\}$ is not H-lsc for $F(x)\cap\{ 1\}=\emptyset$ 
for all $x\in (0,1]$. For some positive results see \cite{Penot-1993}.
\end{enumerate}

Now, we show that if $\mathcal{B}$ is the family of all singletons of $Y$,
or equivalently, the family of all finite subsets of $Y$,
then $\mathcal{B}$-lsc is simply the V-lsc.

\begin{theorem}
Let $\mathcal{B}$ be a cover of $Y$, $F:X\rightarrow 2^Y$,
and consider the following three statements:
\begin{enumerate}
\item $F$ is H-lsc at $x_0$,
\item $F$ is $\mathcal{B}$-lsc at $x_0$,
\item $F$ is V-lsc at $x_0$.
\end{enumerate}
Then $(1)\Rightarrow (2)\Rightarrow (3)$, and the converse implications
does not hold.
\end{theorem}

\begin{proof}
That $(1)$ implies $(2)$ is clear because
$$F(x_0)\cap B\subset F(x_0)\subset W^{-1}(F(x)).$$
Now assume that $F$ is $\mathcal{B}$-lsc at $x_0$. 
In virtue of $\mathcal{B}$-lsc, we may assume that $\mathcal{B}$ contains
all singletons of $Y$. 
We show that $F$ is V-lsc at $x_0$. 
Given an open $G\subset Y$ and $y_0\in F(x_0)\cap G$ we take 
$W\in \mathcal{U}$ such that $W(y_0)\subset G$. 
By the $\mathcal{B}$-lsc there exists a neighbourhood $U(x_0)$ of $x_0$
such that
$$
F(x_0)\cap \{ y_0\}\subset W^{-1}(F(x))\ \text{for all}\ x\in U(x_0).$$
This implies that $y_0\in W^{-1}(y)$ for some $y\in F(x)$,
or equivalently $y\in W(y_0)$ for some $y\in F(x)$.
Consequently, $F(x)\cap G\neq\emptyset$ for all $x\in U(x_0)$.
This proves that $F$ is V-lsc at $x_0$.
It remains to show that the converse implications need not to hold.
Let $Y={\mathbb R}^2$, $\mathcal{B}$ the family of all straight lines of
$Y$ through the origin, and consider the following two
multifunctions:
$F_1(t)$ the line $\{ y=tx\}$, $F_2(t)$ the line $\{ y=1+tx\}$,
for $t\geq 1$.
Observe that:
$F_1$ is V-lsc but not $\mathcal{B}$-lsc at each point,
while
$F_2$ is $\mathcal{B}$-lsc but not H-lsc at each point.
\end{proof}

\begin{theorem}
Let $F:X\rightarrow 2^Y$ be a multifunction and $\mathcal{B}$ be
the family of all singletons of $Y$. Then $F$ is V-lsc at
$x_0$ if and only if it is $\mathcal{B}$-lsc at $x_0$.
\end{theorem}

\begin{proof}
In virtue of the above theorem the second implication is clear.
Now, assume that $F$ is V-lsc at $x_0$ and $\mathcal{B}$ is
the family of all singletons of $Y$. We show that $F$ is
$\mathcal{B}$-lsc at $x_0$. Let $B=\{ y_0\}$. Then $F(x_0)\cap B$
is empty or equal $B$. Let $W\in\mathcal{U}$ be arbitrary. Then
$W(y_0)$ is open and $W(y_0)\cap F(x_0)\neq\emptyset$. By the
V-lsc of $F$ at $x_0$ there exists a neighbourhood $U(x_0)$ of $x_0$ such
that $W(y_0)\cap F(x)\neq\emptyset$ for all $x\in U(x_0)$. Thus for
every $x\in U(x_0)$ there exists $y\in W(y_0)\cap F(x)$. From this
we infer that $(y_0,y)\in W$, or equivalently
$(y,y_0)\in W^{-1}$. Consequently,
$y_0\in W^{-1}(y)\subset W^{-1}(F(x))$. This proves that $F$
is $\mathcal{B}$-lsc at $x_0$.
\end{proof}

A general, similar theorem exists for $\mathcal{B}$-lsc and V-lsc.
First, let us introduce a generalized concept of a totally bounded 
set. Let $\mathcal{B}$ be a cover of $Y$, and $(Y,\mathcal{U})$ 
a quasi-uniform space.
A set $A\subset Y$ will be called $\mathcal{B}$-{\it totally bounded} 
if for every $B\in\mathcal{B}$, the set $A\cap B$ is totally bounded.
Note that, if $\mathcal{B}=\{ Y\}$, then the $\mathcal{B}$-total 
boundedness is simply the usual total boundedness, and if $\mathcal{B}$ 
is the family of all singletons, then the $\mathcal{B}$-total boundedness 
is trivial: each subset of $Y$ is $\mathcal{B}$-totally bounded.

\begin{theorem}
Let $F:X\rightarrow 2^Y$, $x_0\in X$, and $F(x_0)$ be
$\mathcal{B}$-totally bounded, where $\mathcal{B}$ is a cover of $Y$.
Then $F$ is V-lsc at $x_0$ if and only if it is $\mathcal{B}$-lsc at 
$x_0$.
\end{theorem}

\begin{proof}
It is clear that $\mathcal{B}$-lsc implies V-lsc. 
Now, let us assume that $F$ is V-lsc at $x_0$ and $F(x_0)$ is
$\mathcal{B}$-totally bounded. 
Let $B\in\mathcal{B}$ and $W\in\mathcal{U}$ be arbitrary. 
Let $V\in\mathcal{U}$ be such that $V\circ V\subset W$. 
Hence $V^{-1}\circ V^{-1}\subset W^{-1}$, and there exist
$y_1,\ldots,y_n\in F(x_0)$ such that
$$
F(x_0)\cap B \subset 
\bigcup_{i=1}^n V^{-1}(y_i) = 
\bigcup_{i=1}^n V^{-1}(F(x_0)\cap\{ y_i\}).
$$
In virtue of the V-lsc at $x_0$ there exists a neighbourhood $U(x_0)$ of
$x_0$ such that 
$F(x_0)\cap\{ y_i\}\subset V^{-1}(F(x))$, for $i=1,\ldots,n$,
and $x\in U(x_0)$.
This implies that
$F(x_0)\cap B\subset V^{-1}(V^{-1}(F(x))\subset W^{-1}(F(x))$,
which ends the proof.
\end{proof}

As a corollary, for $\mathcal{B}=\{ Y\}$, we get the well-known
equivalence between V-lsc and H-lsc whenever we consider
totally-bounded valued multifunctions.

\subsection*{Remarks}

\subsubsection*{}
One can examine a concept of $\mathcal{B}$-usc described as follows:
for every $B\in \mathcal{B}$ and every $W\in \mathcal{U}$ there
exists a neighbourhood $U(x_0)$ of $x_0$ such that
$F(x)\cap B\subset W(F(x_0))$ for every $x\in U(x_0)$. 
It is clear that V-usc always implies $\mathcal{B}$-usc.
Unfortunately, if $\mathcal{B}$ is the family of all singletons, 
$\mathcal{B}$-usc multifunctions need not to be V-usc.
For example, let $Y={\mathbb R}^2$, $\mathcal{B}$ the family of all 
singletons of $Y$, and we take the following multifunction:
$F(0)$ is the line $\{ y=0\}$, and $F(t)$ is the line $\{ y=tx\}$,
for $t>0$. $F$ is $\mathcal{B}$-usc at 0 but not V-usc at this point.

\subsubsection*{}
See \cite{Costantini-1994} 
for the bounded H-usc, i.e.\ for the $\mathcal{B}$-usc
with $\mathcal{B}$ being the family of all bounded subsets of $Y$.

\subsubsection*{}
It is clear that $\mathcal{B}$-lsc is topologizable 
provided $\mathcal{B}$ is
the family of all singletons or $\mathcal{B}=\{Y\}$.
In general, $\mathcal{B}$-lsc is not topologizable. To show this we
use the following Diagonalization Criterion 
(see e.g.\ \cite{Klein-Thompson-1984}):

\begin{theorem}
Let $T$ be a directed set and for each $t\in T$ there is another
directed set $E(t)$. Then we define a new directed set
$D=T\times\Pi_{t\in T}E(t)$ ordered as follows:
$(t,(\alpha(t)))\leq(s,(\beta(t)))$ if and only if $t\leq s$ and
$\alpha(t)\leq\beta(t)$ for each $t\in T$. Suppose that
$z(t,\gamma)$, $t\in T$, $\gamma\in E(t)$, are elements of
a topological space $Z$. Consider the following net:
$z(t,\alpha)=z(t,\alpha(t))$, $(t,\alpha)\in D$, where
$\alpha(t)$ is the $t$-coordinate of $\alpha$.
If
$\lim_{t\in T}\lim_{\gamma\in E(t)}z(t,\gamma)=z$ then
$\lim_{(t,\alpha)\in D}z(t,\alpha(t))=z$.
\end{theorem}

Now, we can construct an example of a non-topologizable
$\mathcal{B}$-lsc. Observe first that $\mathcal{B}$-lsc is simply
continuity with respect to the following $\mathcal{B}^-$-convergence:
$A_\lambda\rightarrow A_0$ whenever for every $W\in \mathcal{U}$ and
every $B\in\mathcal{B}$ there exists $\lambda_0$ such that
$$
A_0\cap B\subset W^{-1}(A_\lambda)\ \text{for every}\ 
\lambda>\lambda_0.
$$
Let $T={\mathbb N}$ and for each $n\in T$ we take $E(n)={\mathbb N}$.
Let $Y={\mathbb R}^2$ and $\mathcal{B}$ be the family of all straight 
line through the $(0,0)$. For every $n,k\in{\mathbb N}$ we denote:
$A(n,k)=$ the line $\{y=(1/k)x+1/n\}$,
$A(n)=$ the line $\{y=1/n\}$,
$A(0)=B(0)=$ the line $\{y=0\}$.
It is easy to check $\mathcal{B}^-$-convergence:
$\lim_n\lim_kA(n,k)=\lim_nA(n)=A(0)$.

On the other hand the convergence $\lim_{(n,\alpha)}A(n,\alpha(n))=A_0$
does not holds for $A(0)\cap B(0)=$ the line $\{y=0\}$ is not contained
in any $A(n,\alpha(n))+V$, where $V$ is a neighbourhood of $(0,0)$.
Thus the $\mathcal{B}^-$-convergence is not topologizable.

\section{Unions and cartesian products}

In this paragraph we deal with some set-theoretical operations on
multifunctions, namely with unions and cartesian products
(see \cite{Klein-Thompson-1984}).
Operation of intersection of multifunctions will be examined separately
in the next paragraph.

\subsection*{Unions}

Let $X$ and $Y$ be spaces and $F_i:X \rightarrow 2^Y$, $i\in I$,
a family of multifunctions. The {\it union $F=\bigcup_{i\in I}F_i$
of multifunctions} $F_i$ is defined by
$F(x) = \bigcup_{i\in I}F_i(x)$, $x\in X$.
It is known and easy to prove that the union of an arbitrary
family of V-lsc multifunctions is V-lsc. However, the union of an
infinite family of H-lsc multifunctions need not to be H-lsc.
For instance,
define $F_n(0)=[0,n]$ and $F_n(x)=[0,1/ x]$ for $x\in (0,1]$,
$n=1,2,\ldots$. Then the multifunctions $F_n$ are H-lsc at 0.
But, their union $F=\bigcup_n F_n$ is not H-lsc at 0
for $F(0)=[0,+\infty)$.

\begin{theorem}
Let $X$ be a topological space, $(Y,\mathcal{U})$ a quasi-uniform
space and $F_1,F_2:X \rightarrow 2^Y$ multifunctions $\mathcal{B}$-lsc
at $x_0\in X$. Then the union multifunction
$F=F_1\cup F_2$ is $\mathcal{B}$-lsc at $x_0$.
\end{theorem}

\begin{proof}
The proof is a consequence of the formula:
$$W^{-1}(F_1(x))\cup W^{-1}(F_2(x))=W^{-1}(F_1(x)\cup F_2(x)),$$
where $W\in\mathcal{U}$ and $x\in X$.
\end{proof}

It is clear that the above theorem holds for finitely many 
multifunctions, and need not hold when we consider an infinite 
family of multifunctions.

\subsection*{Products}

Now, let us describe the cartesian product of multifunctions.
Let $X$ and $Y_1,Y_2$ be spaces and $F_i:X \rightarrow 2^{Y_i}$,
$i=1,2$, multifunctions. The {\it product} of two multifunctions
$F_1$ and $F_2$ is defined as the multifunction
$F=F_1\times F_2:X\rightarrow 2^{Y_1\times Y_2}$ such that
$F(x)=F_1(x)\times F_2(x)$, $x\in X$.
In particular, if $F_2(x)=Y_2$ for all $x$, or
$F_1(x)=Y_1$ for all $x$, we will write simply, $F_1\times Y_2$,
or $Y_1\times F_1$, respectively.
Analogously, we define the product $\Pi_{i\in I}F_i$ of an arbitrary
family of multifunctions $F_i$, $i\in I$.
It is known that the product of an arbitrary family of V-lsc (H-lsc)
multifunctions is also V-lsc (H-lsc).
Remark that the product of H-lsc multifunctions has more
complicated nature than the product of V-lsc ones.
To formulate a general theorem for $\mathcal{B}$-lsc we need to consider
the product of quasi-uniform spaces.
First, we describe the product of two quasi-uniform spaces.
Let $(Y_i,\mathcal{U}_i)$, $i=1,2$, be quasi-uniform spaces,
and $P_i:Y_1\times Y_2\rightarrow Y_i$, $i=1,2$,
be the projections, i.e.\ $P_i(y_1,y_2)=y_i$, $i=1,2$.
By the {\it product quasi-uniformity
$\mathcal{U}=\mathcal{U}_1\times \mathcal{U}_2$ in $Y=Y_1\times Y_2$ }
we mean a quasi-uniformity generated by the base consisting
of sets of the form
$$
[W_1,W_2]=
\{ (z_1,z_2)\in Y\times Y:(P_i(z_1),P_i(z_2))\in W_i,\ i=1,2 \},
$$
where $W_i\in\mathcal{U}_i$, $i=1,2$. In other words, the set
$[W_1,W_2]$ has a form:
$$
\{ (s_1,y_1,s_2,y_2)\in Y_1\times Y_2\times Y_1\times Y_2:
(s_1,s_2)\in W_1,\, (y_1,y_2)\in W_2\}.
$$
Remark that $[W_1,W_2]=\tilde W_1\cap\tilde W_2$, where
$$
\tilde W_i=
\{ (z_1,z_2)\in Y\times Y:(P_i(z_1),P_i(z_2))\in W_i \},\ i=1,2.
$$
The sets $\tilde W_i$, $W_i\in\mathcal{U}_i$, $i=1,2$, form a subbase
of the product quasi-uniformity $\mathcal{U}_1\times\mathcal{U}_2$.
In case of an arbitrary family of quasi-uniform spaces we
proceed similarly as above and as in the construction of
product topological structures.
Let $(Y_i,\mathcal{U}_i)$, $i\in I$, be a family od quasi-uniform
spaces. Denote: $ Y=\Pi_{i\in I}Y_i$, $P_i$ the projection on
the i-th axis, i.e.\ $P_i(y)=y_i$, where $y_i$ is the i-th
coordinate of $y\in Y$, $i\in I$. By the {\it product quasi-uniformity
$\mathcal{U}=\Pi_{i\in I}\mathcal{U}_i$ in $Y$} we mean a quasi-uniformity
generated by the subbase consisting of sets of the form
$$
\tilde W_i=
\{ (z_1,z_2)\in Y\times Y:(P_i(z_1),P_i(z_2))\in W_i \},
$$
where $W_i\in\mathcal{U}_i$, $i\in I$.
Observe that if $\mathcal{B}_i$ is a cover of $Y_i$, $i\in I$,
then $\Pi_{i\in I}\mathcal{B}_i$, i.e., the family of all sets
of the form $\Pi_{i\in I} B_i$ with $B_i\in\mathcal{B}_i$ and
$B_i=Y_i$ for all but a finite number of $i\in I$, is a cover
of $\Pi_{i\in I}Y_i$.

\begin{theorem}\label{produktowy}
Let $(Y_i,\mathcal{U}_i)$, $i=1,2$, be quasi-uniform spaces
and $A_i\subset Y_i$, $i=1,2$, be arbitrary subsets. Then
\begin{enumerate}
\item
$[W_1,W_2](A_1\times A_2) = W_1(A_1)\times W_2(A_2)$,
\item
$[W_1,W_2]^{-1}(A_1\times A_2) = W_1^{-1}(A_1)\times W_2^{-1}(A_2)$.
\end{enumerate}
\end{theorem}

\begin{proof}
We have
$$
\begin{array}{ll}
\multicolumn{2}{l}{
[W_1,W_2](A_1\times A_2)
}
\\
&
= 
\{ (s_1,s_2)\in Y_1\times Y_2:
\exists\,y_1\in A_1\ \exists\,y_2\in A_2 \,(y_1,y_2,s_1,s_2)\in [W_1,W_2]
\}\\
&
=
\{ (s_1,s_2)\in Y_1\times Y_2: \exists\,y_1\in A_1\ \exists\,y_2\in A_2
\,(y_1,s_1)\in W_1,\,(y_2,s_2)\in W_2 \}\\
&
=
\{ (s_1,s_2)\in Y_1\times Y_2:s_1\in W_1(A_1),\,s_2\in W_2(A_2) \}\\
&
=
W_1(A_1)\times W_2(A_2),
\end{array}
$$
which proves (1). The proof of (2) is similar.
\end{proof}

\begin{theorem}
Let $X$ be a topological space, $(Y_i,\mathcal{U}_i)$, $i\in I$,
a family of quasi-uniform spaces and $F_i:X\rightarrow 2^{Y_i}$
a multifunction $\mathcal{B}_i$-lsc at $x_0\in X$, where
$\mathcal{B}_i$ is a cover of $Y_i$, $i\in I$.
Then the product multifunction $\Pi_{i\in I}F_i$ is
$\Pi_{i\in I}\mathcal{B}_i$-lsc at $x_0$.
\end{theorem}

\begin{proof}
By the construction of the product of quasi-uniform spaces
$(Y_i,\mathcal{U}_i)$, $i\in I$, it is sufficient to consider only
the case $I=\{ 1,2\}$. In general case, the proof is similar.
Let multifunctions $F_i:X\rightarrow 2^{Y_i}$, $i=1,2$, be
$\mathcal{B}_i$-lsc at $x_0\in X$ and $W\in\mathcal{U}_1\times\mathcal{U}_2$
be arbitrary. There exist $W_i\in\mathcal{U}_i$, $i=1,2$, such that
$[W_1,W_2]\subset W$. Now, let $B_i\in \mathcal{B}_1$. There exists
a neighbourhood $U(x_0)$ of $x_0$ such that
$F_i(x_0)\cap B_i\subset W_i^{-1}(F_i(x))$, $i=1,2,\ x\in U(x_0)$,
and, by the Lemma \ref{produktowy}, 
we get
$$
(F_1(x_0)\cap B_1)\times (F_2(x_0)\cap B_2)\subset
[W_1,W_2]^{-1}(F_1(x)\times F_2(x)),
$$
for every $x\in U(x_0)$. 
This shows the $\mathcal{B}_1\times\mathcal{B}_2$-lsc of $F_1\times F_2$
because 
$$
(F_1(x_0)\cap B_1)\times (F_2(x_0)\cap B_2) = 
(F_1(x_0)\times F_2(x_0))\cap (B_1\times B_2).
$$
\end{proof}

\subsection*{Remark}
The converse theorem also holds. Namely, if a product multifunction
$\Pi_{i\in I}F_i$ is $\Pi_{i\in I}\mathcal{B}_i$-lsc at $x_0$, then
for every $i\in I$ the multifunction $F_i$ is $\mathcal{B}_i$-lsc
at $x_0$.

\section{Intersections}

In optimization theory the lower semicontinuity properties of 
intersections of multifunctions play an important role \cite{Penot-1993}.
The most wanted theorems are ones with no boundedness conditions on the
values of intersecting multifunctions. 
Here we formulate a theorem of this kind. 
Let $Y$ be a normed space. 
If we assume that the considered multifunctions are {\it boundedly H-lsc},
i.e., $\mathcal{B}$-lsc with $\mathcal{B}$ being the family of all balls
$B(r)$, $r>0$, then we may formulate a theorem on intersection, without
boundedness conditions on $F(x_0)$.
For this we need a lemma from \cite{Lechicki-Spakowski-1985} on
interiority properties of convex, bounded, and with the nonempty interior
subsets of a normed space.

\begin{theorem}\label{interior-property}
Let $Y$ be a normed space and $A\subset Y$ be convex, bounded, and with
the nonempty interior. 
Then for every $\varepsilon >0$ there exist a set $C\subset int(A)$ and
$\delta >0$ such that $C+B(\delta )\subset A \subset C+B(\varepsilon )$.
\end{theorem}

We need also the following well-known (see \cite{Rabinovich-1967,
Urbanski-1976}) and very useful {\it law of cancellation}:

\begin{theorem}
Let $A$, $B$ and $C$ be subsets of a topological vector space $Y$. 
Assume that $B$ is bounded, and $C$ is nonempty, closed and convex. 
Then 
$A+B\subset \textrm{\it cl\,}(C+B)$ implies $A\subset C$.
In particular, 
$A+B\subset C+B$ implies $A\subset C$, and 
$\textrm{\it cl\,}(A+B)\subset \textrm{\it cl\,}(C+B)$ implies $A\subset C$.
\end{theorem}

Let $X$ be a topological space and $Y$ a topological vector space. 
A multifunction $F$ from $X$ to $Y$ is called {\it locally convex-valued}
({\it locally closed-valued}) {\it at} $x_0\in X$ if there is a 
neighbourhood $U$ of $x_0$ such that $F(x)$ is convex (closed) for every
$x\in U$.

\begin{theorem}
Let $X$ be a topological space, $Y$ a normed space, $\mathcal{B}$ the 
family of balls $B(r)\subset Y$, $0<r$, $F_1$ and $F_2$ two multifunctions
from $X$ to $Y$ and $F=F_1\cap F_2$.
If $F_1$ and $F_2$ are $\mathcal{B}$-lsc at $x_0\in X$, locally convex- 
and locally closed-valued at $x_0$, and 
$\textrm{\it int\,}F(x_0)\neq\emptyset$,
then $F$ is $\mathcal{B}$-lsc at $x_0$, and hence, V-lsc at $x_0$.
\end{theorem}

\begin{proof}
By the assumption on the interior of $F(x_0)$ in $Y$ there
exists $r>0$ such that 
\begin{equation}\label{star0}
\textrm{\it int\,}(F(x_0)\cap B(r))\neq\emptyset.
\end{equation}
Let $\varepsilon >0$ be arbitrary.
By the Lemma \ref{interior-property} there exist a subset 
$C\subset F(x_0)\cap B(r)$ and $\delta >0$ such that
$C+B(\delta )\subset F(x_0)\cap B(r)\subset C+B(\varepsilon )$.
In virtue of the $\mathcal{B}$-lsc at $x_0$ there exists a neighbourhood
$U(x_0)$ of $x_0$ such that 
$$
C+B(\delta )\subset F_i(x_0)\cap B(r)\subset F_i(x)+B(\delta)\
\text{for}\ x\in U(x_0)\ \text{and}\ i=1,2.
$$

We can assume that the multifunctions are closed- and convex-valued
on $U(x_0)$. Applying the law of cancellation, we infer that
$C\subset F_1(x)\cap F_2(x)$ for every $x\in U(x_0)$.
But this implies that 
$$
\begin{array}{lll}
F_1(x_0)\cap F_2(x_0)\cap B(r)&
\subset&
C+B(\varepsilon )\\
&
\subset&
F_1(x)\cap F_2(x)+B(\varepsilon)
\end{array}
$$
for all $x\in U(x_0)$.
This shows that the intersection $F=F_1\cap F_2$ is $\mathcal{B}$-lsc
at $x_0$ and ends the proof.
\end{proof}

\subsection*{Remark}
If $Y$ is finite dimensional we can omit in the above theorem
the assumption that the multifunctions are locally closed-valued,
and then proceed in a manner as in \cite{Lechicki-Spakowski-1985} 
using the below theorem on local interior property of 
$\mathcal{B}$-lower semicontinuous multifunctions.

\section{Vector operations}

Here we consider vector sum and convex hull operations on lower 
semicontinuous multifunctions with values in a topological
vector space (see e.g.\ \cite{Michael-1951}).
Let $X$ be a topological space, $Y$ a topological vector space
and $F,G:X\rightarrow 2^Y$. We define two multifunctions:
$$
\begin{array}{lll}
(F+G)(x)&
=&
F(x)+G(x)\\
&
=&
\{ a+b:a\in F(x),\ b\in G(x) \},\ x\in X,
\end{array}
$$
called the {\it vector sum} of $F$ and $G$, and
$\textrm{\it conv\,}(F)(x)=\textrm{\it conv\,}(F(x))$
the convex hull of $F(x)$, $x\in X$,
called the {\it convex hull} of $F$.

\subsection*{Vector sum}
It is known and easy to prove that the vector sum of two H-lsc
multifunctions is H-lsc. We state some further results and show that,
in general, the vector sum of two $\mathcal{B}$-lsc multifunctions 
need not to be $\mathcal{B}$-lsc.

\begin{theorem}
Let $F,G:X\rightarrow 2^Y$ be V-lsc at $x_0\in X$.
Then the vector sum $F+G$ is V-lsc at $x_0$.
\end{theorem}

\begin{proof}
Let $V$ be an arbitrary neighbourhood of 0 in $Y$, $F$ and $G$ be V-lsc at
$x_0$, and recall that V-lsc is equivalent to $\mathcal{B}$-lsc with 
$\mathcal{B}$ equals the family of all singletons of $Y$. 
Let $b\in Y$ be such that $(F(x_0)+G(x_0))\cap\{ b\}\neq\emptyset$. 
Then $b=b_1+b_2$ with $b_1\in F(x_0)$ and $b_2\in G(x_0)$. 
By the $\mathcal{B}$-lsc of $F$ and $G$ there exists a neighbourhood
$U(x_0)$ such that
$$
\begin{array}{l}
F(x_0)\cap\{ b_1\}\subset F(x)+V\ \text{and}\\
G(x_0)\cap\{ b_2\}\subset G(x)+V
\end{array}
$$
for all $x\in U(x_0)$. 
This implies that 
$$b_1+b_2\in F(x)+G(x)+V+V$$ 
for $x\in U(x_0)$,
or equivalently,
$$(F(x_0)+G(x_0))\cap\{ b\}\subset F(x)+G(x)+V+V$$ 
for $x\in U(x_0)$,
which shows $\mathcal{B}$-lsc of $F+G$ at $x_0$ and ends the proof.
\end{proof}

We say that a cover $\mathcal{B}$ of $Y$ is {\it translation invariant}
if $B+c\in\mathcal{B}$ for every $B\in\mathcal{B}$ and every vector $c\in Y$.

\begin{theorem}
Let $\mathcal{B}$ be a translation invariant cover of $Y$,
$F:X\rightarrow 2^Y$ a multifunction $\mathcal{B}$-lsc at $x_0\in X$
and $g:X\rightarrow Y$ a function continuous at $x_0$. Then
the vector sum $(F+g)(x)=F(x)+g(x)$, $x\in X$, is a multifunction
$\mathcal{B}$-lsc at $x_0$.
\end{theorem}

\begin{proof}
Let $V$ be a neighbourhood of 0 in $Y$ and $B\in\mathcal{B}$. 
Note that 
$$(F(x_0)+g(x_0))\cap B=F(x_0)\cap(B-g(x_0))+g(x_0)$$
By the assumptions for all $x$ in a neighbourhood of $x_0$ we have
$$\begin{array}{lll}
(F(x_0)+g(x_0))\cap B&
\subset&
(F(x)+g(x_0)+V)\cap B\\
&
\subset&
(F(x)+g(x)+V+V).
\end{array}$$
This shows $\mathcal{B}$-lsc of $F+g$ at $x_0$ and ends the proof.
\end{proof}

The following example shows that translation invariantness of
$\mathcal{B}$ is not sufficient to get $\mathcal{B}$-lsc of the vector
sum of a two $\mathcal{B}$-lsc multifunctions.

\begin{example}
Let $Y={\mathbb R}^3$, $\mathcal{B}$ consists only of the plane $\{ y=x\}$
and all of its translations. Define two multifunctions:
$F(t)$ the line $\{ z=ty,\, x=0\}$ and
$G(t)$ the line $\{ y=0,\, z=0\}$, $t\geq 0$.
Observe that $F$ and $G$ are $\mathcal{B}$-lsc but their vector sum
$F+G$ is not $\mathcal{B}$-lsc at each point.
\end{example}

\subsection*{Convex hull} 
It is known \cite{Michael-1956} and easy to proof 
that for every V-lsc multifunction $F$, the convex hull of $F$ is also 
V-lsc. We use the concept of $\mathcal{B}$-lsc to get a general result
provided the space $Y$ is locally convex. In particular, we get a result
for H-lsc.

\begin{theorem}
Let $Y$ be a locally convex space and $F:X\rightarrow 2^Y$ be
$\mathcal{B}$-lsc at $x_0\in X$.
Then the convex hull of $F$ is $\mathcal{B}$-lsc at $x_0$.
\end{theorem}

\begin{proof}
Let $V\subset Y$ be a convex neighbourhood of 0 in $Y$ and
$B\in\mathcal{B}$. By the $\mathcal{B}$-lsc of $F$ at $x_0$ there exists
a neighbourhood $U(x_0)$ of $x_0$ such that
\begin{equation}\label{star}
F(x_0)\cap B\subset F(x)+V\ \text{for all}\ x\in U(x_0).
\end{equation}
We claim that 
$$
\textit{conv\,}(F(x_0))\cap B\subset\textit{conv\,}(F(x))+V\
\text{for all}\ x\in U(x_0).
$$ 
Indeed, let $y\in\textit{conv\,}(F(x_0))\cap B$ and $x\in U(x_0)$ be
arbitrary. 
In virtue of (\ref{star}) there exist 
$n\in{\mathbb N}$,
$y_1,\ldots,y_n\in F(x_0)$, 
$z_1,\ldots,z_n\in F(x)$,
$v_1,\ldots,v_n\in V$, and positive numbers 
$t_1,\ldots,t_n$
such that 
$y=t_1y_1+\cdots +t_ny_n$, 
$t_1+\cdots +t_n=1$, and 
$y_i=z_i+v_i$ for every $1\leq i\leq n$. 
This implies that
$$
y=t_1z_1+\cdots +t_nz_n+t_1v_1+\cdots +t_nv_n\in \textit{conv\,}(F(x))+V
$$
for $V$ is convex, which ends the proof.
\end{proof}

\subsection*{Remark}
If the topological vector space $Y$ is not locally convex then the
convex hull operation does not preserve H-lsc. Indeed, if $Y$ is 
metrizable and not locally convex then there exists a sequence
$y_n\in Y$ which converges to $0$ such that the convex hull of the
set $\{y_n:n=1,2,\ldots\}$ is not bounded \cite{Araki-1995}. 
Now, observe that the multifunction $F$ defined by: 
$F(0)=\{0,y_1,y_2,\ldots\}$, 
$F(1/n)=\{0,y_1,\ldots,y_n\}$, $n=1,2,\ldots$, is H-lsc at $0$ but 
the convex hull of $F$ is not.
For a simple example of such sequence $y_n\in l^p$ ($0<p<1$) see 
\cite{Rudin-1973} or \cite{Araki-1995}.

\providecommand{\bysame}{\leavevmode\hbox to3em{\hrulefill}\thinspace}
\providecommand{\MR}{\relax\ifhmode\unskip\space\fi MR }
\providecommand{\MRhref}[2]{%
  \href{http://www.ams.org/mathscinet-getitem?mr=#1}{#2}
}
\providecommand{\href}[2]{#2}

\end{document}